# Paramathematical notions and Klein's Plan B: the case of equations

Carl Winsløw[1]

[1]University of Copenhagen, Faculty of Science, Copenhagen, Denmark; winslow@ind.ku.dk

*Undergraduate students of mathematics continue to solve equations in virtually any course they attend, just as they did in secondary school – yet what do they learn about equations and their solutions at university? Are they capable to combine elements of abstract algebra and real analysis to assess what it means to solve an equation, in particular, what it means for an equation to be solvable – or not? In this paper, we present a mathematical and didactical analysis of these questions, illustrated by examples from a capstone course for future Danish high school teachers.*

*Keywords: Equations, Klein's Plan B, Analysis, Algebra, paramathematical notions*

## Introduction: what is an equation?

While every school child has some idea of what equations are – based on examples encountered from primary school and onwards – it is harder both for them and many others to provide a formal mathematical definition. According to Wikipedia (n.d., underlining is by the author):

> In mathematics, an equation is a <u>mathematical formula</u> that expresses the equality of two <u>expressions</u>, by connecting them with the equals sign =. (…) Solving an equation containing <u>variables</u> consists of determining which <u>values</u> of the variables make the <u>equality</u> true.

The underlined terms are more or less undefined; for instance, the above text links to an explanation of "mathematical formula" which refers back to equations, and provides the example

$$V = \frac{4}{3}\pi r^3 \quad (1)$$

for the volume of a sphere. This is an example of a "formula" which is also a "theorem", given that it is true with an appropriate interpretation of the terms. However, for this type of equation, it is not immediately clear what "solving" means, since – with the said interpretation – it is *always* true. In the above quotation, the other terms are equally difficult to define in a precise, non-circular manner. The undergraduate student will encounter many other kinds of "formulae" connecting "expressions" by one or more equals signs, and the values obtained by "solving" such "equations" are no longer limited to (real) numbers, but could be group elements, vectors, functions and so on. While students are not presented with an explicit definition of the underlined terms – or of the notion of equation itself – the terms are nevertheless extensively used both in school and at university. They are what Chevallard (1985, pp. 49-51) calls *paramathematical notions* : they are named and used as tools in the study of explicitly defined mathematical objects, but they are not themselves explicitly defined.

Among the more commonly encountered types of equations in middle and high school, we can of course point to familiar cases such as linear, quadratic or (more generally) *polynomial* equations. These can be defined as expressions of type $p(x) = 0$, where the left-hand side is a polynomial over some set $E$ equipped with two operations satisfying, for instance, the axioms of an integral domain (the cases of integers, rational or real numbers with the usual operations being the most familiar ones). Here, the notion of solution takes on the more explicit meaning $\{x \in E : p(x) = 0\}$, the subset of $E$

(or the "values") that makes the equation true. A more general definition of what is most often thought of as an equation in secondary level mathematics can be given as follows (in terms of genuine mathematical notions): given a subset $E$ of $\mathbb{R}$ and a function $f: E \to \mathbb{R}$, solving "$f(x) = 0$" means to determine the subset $\{x \in E : f(x) = 0\}$ of $E$. Here, we can subsume the variation "$f(x) = g(x)$", with another function $g: E \to \mathbb{R}$ on the right-hand side, since $f - g$ is also a function. The variation – with more than one variable, say, $x_1, \ldots, x_n$, and a product set $E_1 \times \ldots \times E_n$ in which solutions are to be found – is not without relevance even in school (cf. (1), if we consider $V$ and $r$ as mere variables) and does not add any substantial complexity to task of *defining* equations and the meaning of solution in this more restrictive real-valued function context. Thus, we shall stay here with the above, simpler definition of what we could call "a functional equation in one real variable with domain $E \subseteq \mathbb{R}$".

University students are certainly accustomed to solving such equations by algebraic manipulation. What these "manipulation" are, or could be, is much less clear – and even less so, what conditions allow them to complete the solution. We now delve into observations from a practical context.

## Research questions and context of observations

In many educational systems, secondary level teachers are mainly prepared at universities through standard university courses in their subjects. Such courses are typically more designed to prepare for a scholarly career than for teaching at secondary level, and mathematics may be one of the subjects where the distance between even introductory university courses (on subjects like real analysis and abstract algebra) and secondary level mathematics is the greatest. Already in 1908, Klein (2016) pointed to this "discontinuity" and emphasized the responsibility of universities to help future teachers connect their university studies to the mathematics they will be teaching. In a series of papers (from Winsløw, 2013 to Winsløw and Huo, to appear), we have investigated this from the point of view of *institutional relationship to knowledge* (Chevallard, 1985), namely as the development a new (teaching oriented) relationship to the mathematics of secondary school, at the end of university mathematics studies. As in these previous papers, the context from which the examples in this paper are drawn comes from the course "Mathematics in a teaching context" (abbreviated UvMat from the Danish name) which students taking a minor (second subject) degree in mathematics attend towards the end of their studies at the University of Copenhagen. In this paper, we focus on *how elements of abstract algebra and real analysis may be drawn upon to supply students with a new and more advanced relationship to the notion of equation*. Courses for future teachers seem an obvious occasion to take up such paramathematical notions, but the scope of our hypothesis is not limited to future teachers. We propose that paramathematical notions could be usefully revisited at the end of undergraduate mathematics studies in general, in a setting which Durel (1993, 223) explains thus:

> The capstone course typically is defined as a crowning course or experience coming at the end of a sequence of courses with the specific objective of integrating a body of relatively fragmented knowledge into a unified whole. As a rite of passage, this course provides an experience through which undergraduate students both look back over their undergraduate curriculum in an effort to make sense of that experience, and look forward to a life by building on that experience.

The idea of *integrating a body of relatively fragmented knowledge* rejoins another central idea of Klein (2016, 77-85), namely that of two different "plans" for creating and learning mathematics (cf.

also Kondratieva & Winsløw, 2018): the first, Plan A, consists in building up highly structured knowledge in separate domains, while the second, Plan B, consists in combining and integrating knowledge from different domains. We shall argue, through the example of equations, that advancing students' knowledge of central paramathematical notions can benefit from a Plan B approach. Given the preponderance of "solving equations" in the teaching of mathematics at pre-university levels, it is naturally one which teachers may particularly benefit from getting to appreciate more clearly, through a well-designed capstone experience. Also for the future scholar, it could be useful to acquire a more integrated view of this and other paramathematical notions, like variable, structure and proof.

In UvMat, students' preparation amounts to at least 1½ years of undergraduate courses, in subjects such as discrete mathematics, probability, analysis and abstract algebra. As some of the observations we shall share here suggest, they typically do not have full command of the formal content of these courses; most of them do their major in a subject from the humanities and will become teachers. An important function of UvMat is therefore to revisit and combine selected elements of the prerequisite courses in order to consolidate and test their knowledge of fundamental mathematical – and paramathematical – notions that appear directly in the secondary curriculum. As in other courses, challenging tasks form an important component of the course, along with theory that integrates and develop their prerequisite knowledge (see Winsløw and Huo, to appear, concerning the mathematical object $\mathbb{R}$). In particular, every week students get a mandatory "weekly assignment", which develops, through six questions, certain theoretical complements to the mathematical contents that have been studied in the week's lectures and exercise sessions. The course is assessed through a written exam.

## A "tricky" equation and students' struggle with it

The most basic way to approach the topic of equations is naturally to discuss concrete and simple cases, with tasks asking students to solve an equation while explaining each step (technique) applied in the solution. At university level, going beyond standard secondary level types of tasks (like quadratic equations) is both possible and expected, for instance to help students discover that certain standard techniques may need caution in some cases. As an example, the following exercise was one among five proposed at the final five-hour exam of the 2023 edition of UvMat:

> Let the function $f$ be given by $f(x) = \sqrt{x} + \sqrt{2x+1}$ for $x \geq 0$.
>
> a) Show that if $a \geq 1$, the equation $f(x) = a$ has precisely one solution. (Hint: is $f$ increasing?)
>
> b) One of your students solved the equation $f(x) = 3$ manually, and got two solutions. How could that be?

In a), students are expected to show that $f$ is in fact strictly increasing (being a sum of two strictly increasing functions) and therefore injective; so there is at most one solution. Using the continuity of $f$ and the fact that $f(x) \to \infty$ as $x \to \infty$, one further shows that $f([0, \infty[) = [1, \infty[$, and existence follows. To solve the equation in b), one can successively manipulate the equation to get $x + 2x + 1 + 2\sqrt{x(2x+1)} = 9$ (squaring), $3x - 8 = -2\sqrt{2x^2 + x}$ (reordering), and then get $(3x-8)^2 = 4(2x^2 + x)$ (squaring), which is a quadratic equation with two solutions. However, the last equation is not equivalent to second one, which implies $3x - 8 \leq 0$, while only one of the two solutions to the last equation satisfies this. In other words, the squared equation is implied by, but not equivalent to, the given equation. Thus, the imaginary student may have used correct algebraic techniques, but forgotten a point which is (at the level of theory) emphasized in the course: when solving an equation by step-by-step rewritings, logical equivalence must be attended to, to avoid losing solutions or, as is in this case, adding false ones. Implication corresponds in fact to set inclusion: here, $f(x) = 3 \Rightarrow (3x - 8)^2 = 4(2x^2 + x)$ means that $\{x \geq 0: f(x) = 3\} \subseteq \{x \geq 0: (3x - 8)^2 = 4(2x^2 + x)\}$, but the reverse implication and inclusion do not hold.

At the exam, only about half of the students provided a complete argument for a), while less than half were able to give a complete and correct answer to b). As for a), many students overlooked the need to prove existence of solutions (and not merely uniqueness, which is supported by the hint); thus, the contribution of continuity (from analysis) was not activated, despite several similar examples having been treated in the course. In b), most students carry out the rewritings indicated above, to arrive at two solutions, but fail to explain why these are not both correct (except by reference to the result given in a)). A few students did not even master the algebraic techniques, as one student who writes:

> If the student had solved the equation correctly, it would look somewhat like this [*see Figure 1 for the calculations*]. And the student would then have seen that there is only one solution to this equation. But the student possibly made some calculation error, and arrived at a quadratic equation, from which two solutions are obtained.

Figure 1: Student rewriting attempt

The student apparently uses the "rule" $(a + b)^2 = a^2 + b^2$ (or, perhaps, $\sqrt{a + b} = \sqrt{a} + \sqrt{b}$) and does not explain further how "calculation error" could lead to a quadratic equation. Most likely, the student did not check the solution by insertion into the original equation, despite having 1 hour to answer each exercise. It is not surprising that students, who made such errors, performed poorly on other tasks as well, and typically failed the exam. What is more surprising is that they could pass several undergraduate courses as explained above, and still produce such solutions. Indeed, a capstone course may not only consolidate and enrich fundamental knowledge, but also be a final and necessary test of it. In fact, the solution of b) does not draw significantly on university courses in algebra and analysis.

### Algebra and analysis in the undergraduate curriculum

As documented by Bosch et al. (2021), undergraduate mathematics programmes in Europe share many common traits, including two dominant sequences of courses: one in linear and abstract algebra, and one in calculus and analysis. While linear algebra is heavily drawn on in the analysis sequence, advanced analysis and abstract algebra appear relatively independent at this level (while at the graduate level, links do appear, like applications of group and ring theory in more advanced analysis courses). As one German mathematician cited in the above paper (p. 148) said,

> You do obtain a solid background in analysis and algebra [and] it gives a solid foundation for most of the mathematical directions…

But these "directions" are not really reached by those who take only the foundational (undergraduate) courses. The polynomials encountered in algebra (for instance as elements of key examples of rings, or in Galois theory) are not necessarily connected, by the students, to the polynomial functions which appear here and there in real and complex analysis (for instance in Weierstrass' approximation theorem, or in the "fundamental theorem of algebra", which is not really a theorem of algebra, and is usually proved in a first course on complex analysis).

The situation is, as the above suggests, not entirely symmetric. Abstract algebra is often presented as a sequence of derivations within closed axiomatic theories; analysis is not drawn on in this development, but may be referred to when examples related to the real and complex number fields are considered. The theoretical foundations of analysis involve both algebraic and topological structures that can, in principle, be developed similarly (see Hausberger and Hochmuth, to appear). But at the undergraduate level, these foundations are usually presented in more informal ways, especially when it comes to the fundamental properties of the real numbers (Winsløw and Huo, to appear). With this asymmetry, undergraduate presentations of abstract algebra, and of post-calculus analysis, appear as more or less independent "paths". For future teachers in particular, the analysis path formalizes and extends crucial parts of the calculus taught in secondary schools, while abstract algebra appears mostly as a new, self-contained and much less familiar universe of mathematical results. As we shall now see, the paramathematical notion of equations can help complementing this Plan A approach with elements of Plan B (both in the sense of Klein, explained above).

## What are "unsolvable" equations?

In abstract algebra courses, students encounter the notion of algebraic number, defined as roots of non-zero polynomials with integer coefficients. While these numbers are all complex numbers, the construct of field extension means that the algebraic numbers can be studied with algebraic means, without a need of *all* complex numbers; and it remains a result of abstract algebra that a polynomial $p$ is solvable by radicals (loosely, its roots can be expressed using only rational expressions of its coefficients and of radicals of elements of the underlying field) if and only if its Galois group is solvable. This leads to the impossibility of generalising algebraic methods (as for solving polynomial equations of degree 2) beyond degree 4. While this result is rarely proved at the undergraduate level, it if often stated and partially developed there, to show that even polynomial equations are, in general, significantly more complicated than the special cases known from secondary school.

The notion of algebraic number is naturally worth elaborating on. In UvMat, students review countable and uncountable sets (also encountered in analysis courses). Based on a semi-formal construction of $\mathbb{R}$ as the set of infinite decimals (with some, limited, identifications), they are exposed to a semi-formal proof that $\mathbb{R}$ is uncountable. They also solve an exercise set up to prove that the set of algebraic, real numbers is countable. As a result, they know that there are plenty of non-algebraic (or *transcendental*) real numbers. In some years, they are exposed to a weekly assignment which leads them through a proof that Euler's number $e$ is indeed transcendental. The proof relies only on methods from first year courses in analysis, and a few basic facts about prime numbers.

Coming back to equations and functions, the definition of algebraic numbers can be generalised to algebraic *functions*: a function $f: D \to \mathbb{R}$ is said to be *algebraic* if there is a non-zero polynomial $p$ of two variables, with integer coefficients, such that $p(x, f(x)) = 0$ for all $x \in D$. A function is called *transcendental* if it is not algebraic. These definitions are, of course, purely algebraic.

In a weekly assignment, these definitions were provided and the first question for students was:

a) Prove that $f: D \to \mathbb{R}$ is algebraic if and only if there exist (one-variable) polynomials $p_0,\ldots,p_n$ with integer coefficients, such that not all $p_i$ are zero, and such that $\sum_{k=0}^{n} p_k(x) f(x)^k = 0$ for all $x \in D$.

Question b) then asked to explain that some basic examples, like polynomials over $\mathbb{Z}$ and the function $f(x) = \sqrt[3]{2}x + \sqrt{5}$, are algebraic. Both a) and b) turned out to be very challenging for most students, and almost all had to redo this part of the assignment (as it is a requirement for approval that all questions are answered correctly). Question c) asked to show that if $a$ is a transcendental real number, then $f(x) = a^x$ defines a transcendental function; question d), that if $f$ is algebraic and injective, then so is $f^{-1}$. Here, a main challenge seemed to be to choose a convenient characterization to work with. Using the definition, c) is easy: roughly, $p(x, f(x)) = 0$ if and only if $p(f^{-1}(x), x) = 0$. To answer d), the characterisation from a) comes in handy since $f(1) = a$. As an application of c) and d), the students easily show that the natural logarithm is a transcendental function (question e)).

Finally, students were presented with the standard definition of *transcendental equation*, being of the form $f(x) = g(x)$ where at least one of $f$ and $g$ is a transcendental function, and the last question f) asks them to give two examples of transcendental equations with solutions, one which they can solve, and one which they cannot. They were also asked to solve both equations numerically (with a computer algebra system). Given that students had just seen examples of both algebraic and transcendental functions, question f) was not difficult: for instance, $e^x = \frac{1}{2}$ and $e^x = x + 2$ work.

The last question was further discussed in a follow-up session, as some students were rightly intrigued that transcendental equations are sometimes easy or even trivial (like $e^x = e^x$), while algebraic equations (polynomial, for instance) may indeed be impossible to solve algebraically. I now outline some of the ideas this led to develop, in part with the students, in part while elaborating this paper.

Of course, "difficulty" of a task is not even a paramathematical notion, but depends on the knowledge of the solver and context – more precisely, the means available and permitted for solving. In the case $e^x = x + 2$ mentioned above, a professional tool like Maple does give two explicit (real) solutions:

*solve*(exp(*x*) = *x* + 2, *real*);

-LambertW(-exp(-2)) - 2, -LambertW(-1, -exp(-2)) - 2

The only problem with this solution is that most undergraduate students have never encountered the functions which Maple calls "LambertW" and "LambertW(−1)", while the rest makes perfectly sense to them. Indeed, university students' inventory of (special) functions defined on $\mathbb{R}$ is not extended much, if at all, beyond what they already encountered in secondary school. The definition of algebraic function merely points to roots of polynomial (functional) equations, many of which could indeed be quite foreign to students (an example is given below), and the realm of transcendental functions is even vaster. Introducing LambertW or at least its principal branch, as has been done in some weekly assignments for other years, does not solve the general problem, but suggests that there is a certain arbitrarity about the inventory of functions that become "familiar" in secondary and undergraduate mathematics. We now reconsider the relation between functions and equations to explain how this inventory affects their conception of what it means to (be able to) solve an equation.

## Solving equations from algebraic, analytic and didactic perspectives

Let us return to our definition of equations in one real variable, and what it means to solve them: an equation can be written on the form $f(x) = 0$ where $f: D \to \mathbb{R}$ is a function and $D \subseteq \mathbb{R}$; solving such an equation amount to determine $N_f = \{x \in D: f(x) = 0\}$. What is now at stake is what are permissible ways to determine $N_f$. It is clear that simply stating the set $N_f$ will not qualify as a valid solution; $N_f$ must be given in a somehow explicit form, along with an argument of its validity. In the algebraic realm, where only the algebraic (and not the topological) structure of $\mathbb{R}$ is considered, we cannot use properties based on completeness, such as the implicit function theorem. Functions defined only by means of the algebraic structure are all rational; if $f = p/q$ where $p$ and $q$ are non-zero polynomials, we have $N_f = N_p \backslash N_q$. Thus, the solution of algebraic equations amounts to finding the roots of two polynomials. Doing so by algebraic means amounts to specifying the (finitely many) roots of $p$ and $q$ by carrying out explicit calculations in the splitting fields, resulting in explicit formulae of the roots in terms of the coefficients of the polynomials (here, roots and rational expressions are allowed). As we have noted, this cannot be done for all polynomials of degree higher than four – for instance, it cannot be done for $p(x) = x^5 - x - 1$. Still, there are many common polynomial equations of degree higher than four that can be solved using "algebraic manipulation" (for a somewhat incomplete attempt to formalize this paramathematical notion, see Weiss, 2020, chap. 3). For instance, $x^6 - x^3 - 1 = 0$ can be solved using the substitution $y = x^3$.

Staying with polynomial functions and equations is clearly inadequate for both future mathematicians and teachers. If in our equations, we allow functions defined by analytic means – such as $f(x) = e^x$ – then one may not reasonably expect algebraic definitions (like that of algebraic equation) to provide a satisfying answer to the meaning of solvable equation, of acceptable specification of the solution, or of acceptable proof of it. A simple and explicit way to define $f(x) = e^x$ is as the inverse function to $g(x) = \int_1^x \frac{dt}{t}$ (defined for $x > 0$). This involves proving that $g$ is injective and has range $\mathbb{R}$, so that the inverse exists and is defined on $\mathbb{R}$. Then, the equation $e^x = \frac{1}{2}$ has the perfectly explicit solution $\{g(\frac{1}{2})\}$. But once such functions are allowed – as they are even in secondary school! – it must be just as admissible to define the inverses of the function $F(x) = x \cdot e^x$ corresponding to the two intervals where $F$ is injective. Then one can easily get the solutions to $e^x = x + 2$ that were obtained by Maple above, by observing that the equation is equivalent to $-e^{-2} = -(x+2)e^{-(x+2)}$.

Allowing such arguments also ends all difficulties with polynomial equations. For instance, basic analysis shows that with $p(x) = x^5 - x - 1$, there is a unique $\alpha > 1$ such that $p(\alpha) = 0$. What remains is to use the injectivity of the restriction of $p$ to $]1, \infty[$ to define its inverse $Z$ on $]-1, \infty[$ and conclude $N_p = \{Z(0)\}$. From an analytic viewpoint, nothing makes $Z$ a less respectable function than the natural logarithm. From an algebraic perspective, $Z$ is in fact algebraic (as is the equation just solved).

From a *didactic* perspective, it remains that students have been systematically trained to consider a small set of transcendental functions (like exp, ln, sin, cos a few more) as "concrete" functions, even if in many cases they have not even been given a concrete and complete definition of these. These functions, in combinations with algebraic functions and a few transcendental numbers like $e$ and $\pi$, constitute all the concrete functions they have ever encountered. As a result, $\sin \frac{\sqrt[3]{x+1}}{x^2 + \pi}$ is a concrete function to them, while *Z* and *LambertW* are not. In this sense, after their exposure to real analysis, it seems useful to both (re)study the definitions of the functions they are familiar with (cf. Winsløw, 2013) and to work with less familiar – but often widely applicable – new functions. Both efforts will extend their mathematical notion of function and, what is the matter considered here, provide a firmer basis to their paramathematical notions of equation and what it means to solve one. For the latter purpose, the distinction between algebraic and transcendental equations is at most the beginning.

## Conclusion

We have outlined, though an example, how a capstone course can *question* and *combine* knowledge from major areas of the undergraduate study, in view of strengthening the students' autonomous perspective on elementary and everywhere present mathematical and paramathematical notions. The case of functions and equations shows that students may not only have gaps in their knowledge from major mathematical domains like algebra and analysis, but that even the knowledge they have can be quite compartmentalized. It may include theoretical results that have highly underdeveloped links to more elementary and foundational notions, like the ones considered here. Klein's Plan B therefore has all its place as a central didactical challenge for capstone courses in mathematics.